\documentclass[12pt]{amsart}

\usepackage[notref,notcite]{showkeys}

\newtheorem{thm}{Theorem}[section]
\newtheorem{prop}[thm]{Proposition}
\newtheorem{defin}[thm]{Definition}

\newtheorem{corol}[thm]{Corollary}

\newtheorem{lem}[thm]{Lemma}
\newtheorem{conj}[thm]{Conjecture}
\newtheorem{rem}[thm]{Remark}

\newcommand{\N}{{\mathbb N}}

\newcommand{\pp}{{\mathcal P}}

\newcommand{\R}{{\mathbb R}}\newcommand{\RR}{{\mathbb R}^2}

\newcommand{\ep}{\varepsilon}

\newcommand{\al}{\alpha}
\newcommand{\be}{\beta}
\newcommand{\ga}{\gamma}\newcommand{\Ga}{\Gamma}
\newcommand{\de}{\delta}

\newcommand{\la}{\lambda}

\begin{document}

\bibliographystyle{plain}

\title[Entropy and insecurity]
{Growth of the number of geodesics between points and insecurity for riemannian manifolds}

\author{Keith Burns and Eugene Gutkin}

\address{Department of Mathematics\\
Northwestern University\\
Evanston, IL 60208\\
USA}

\email{burns@math.northwestern.edu}
 
\address{IMPA \\
Estrada Dona Castorina 110 \\
Rio de Janeiro,   22460-320\\
Brasil}

\email{gutkin@impa.br}

\thanks{K.B.~is partially supported by N.S.F.~grant DMS-0408704.}

\date{\today}

\begin{abstract} 
A Riemannian manifold is said to be uniformly secure if there  is a finite
number $s$ such that all geodesics connecting an arbitrary
pair of points in the manifold can be blocked by $s$ point obstacles.
We prove that the number of geodesics with length $\leq T$ between
every pair of points in a uniformly secure manifold grows
polynomially as $T \to \infty$. By results of Gromov and Ma\~n\'e, 
the fundamental group of such a manifold  is virtually nilpotent, and the topological 
entropy of its geodesic flow is zero. We derive from this that a
compact Riemannian manifold with no conjugate points whose geodesic flow has 
positive topological entropy is totally insecure: the geodesics between any 
pair of points cannot be blocked by a finite number of point obstacles.   
\end{abstract}

\maketitle

\section{Introduction}       \label{intro}

This paper explores interactions between two approaches to the study of geodesics 
in a Riemannian manifold. Both approaches study geodesic segments
between points of the manifold. One of them relates
the growth as $T \to \infty$
of the number of these geodesics with length $\leq T$ to the topological
entropy of the geodesic flow. The other studies whether
or not these geodesics can be blocked by a finite number of point obstacles.

Let $M$ be a compact Riemannian manifold  that is uniformly 
secure with security threshold $s$. 
This means that the geodesics between any two 
points in $M$ can be blocked by $s$ obstacles.
Our main technical result concerns the growth rate of joining geodesics
in such a manifold. We show that there are constants
$C$ and $d$, which depend only on $s$ and the injectivity radius, 
such that the number of geodesics 
with length $\leq T$ between any pair of points in $M$ is bounded from above
by $CT^d$. See Proposition~\ref{main_prop} and Theorem~\ref{polin_thm}. 
This bound has strong consequences for the geometry and 
topology of $M$. In particular, its fundamental group is virtually
nilpotent. See Theorem~\ref{entrop_thm}.


Our methods yield strong results for manifolds with no conjugate points.
In particular, a compact manifold with nonpositive curvature
satisfies the following dichotomy.
Either such a manifold is flat and the geodesics between
every pair of points are blocked by a finite number of obstacles; 
or the manifold is not flat and removing a finite number of points from the manifold
cannot block the geodesics connecting any pair of points. See Corollary~\ref{neg_curv_cor}.

\medskip

We thank J.F.\ Lafont and B.\ Schmidt for pointing out a simplification 
to our original proof of Theorem~\ref{nocp_thm}. They have independently
studied connections between security and the geometry of a manifold \cite{LS}.
In particular, they have obtained Theorem~\ref{nocp_thm} and  the lemmas in 
Section~\ref{keysect}.

\section{Background and Previous Results} \label{back&prev}
Throughout this paper by a (Riemannian) manifold we will mean a 
complete, connected, infinitely differentiable (Riemannian) manifold. 
Geodesics will always be parametrized by arclength. Unless otherwise specified, 
a geodesic will be a segment of finite, positive length.
Let $M$ be a  Riemannian manifold.
If $[a,b]$ is a closed bounded interval and $\ga: [a,b] \to M$ is a geodesic, the 
{\em endpoints} of $\ga$ are its {\em initial point} $\ga(a)$ and its {\em final point} $\ga(b)$.
The points $\ga(s)$,
where $a < s < b$, are the {\em interior points} of $\ga$. 
We will say {\em $\ga$ passes through $z$} to indicate that 
$z$ is an interior point of $\ga$.
Note that a point $z\in M$ may be both an interior point and an endpoint of $\ga$.

A {\em configuration} in $M$ is an ordered pair
of points in $M$. The points in question may coincide.
\begin{defin}      \label{connect_def} 
{\rm Let $(x,y)$  be a configuration in $M$. 
A geodesic $\ga$ {\em joins} $x$ to $y$ if $x$ is its initial point and
$y$ is its final point. 
A geodesic $\ga$  {\em connects} $x$ and $y$ if it joins $x$ and $y$
and does not pass through either $x$ or $y$.}
\end{defin}

We denote the set of geodesics joining $x$ and $y$ by $G(x,y)$
and the set of geodesics connecting $x$ and $y$ by $\Gamma(x,y)$.
Of course $\Gamma(x,y) \subset G(x,y)$. The subsets of $G(x,y)$ and
$\Gamma(x,y)$ consisting of geodesics with length $\leq T$ will be denoted
by $G_T(x,y)$ and $\Gamma_T(x,y)$ respectively. 
Let $m_T(x,y)$ be the number of geodesics in $\Gamma_T(x,y)$ 
and let $n_T(x,y)$ be the number of geodesics in $G_T(x,y)$.
If one of the sets is infinite, the corresponding number will be $\infty$.

\subsection{Growth of $n_T(x,y)$}                 \label{SSn_T} 
The behavior  of 
$n_T(x,y)$ for compact  manifolds has been extensively 
studied beginning with Berger and Bott's 1962 paper \cite{BeBo}. 
They observed that for each $T > 0$ the function $n_T(x,y)$ is 
finite and locally constant on an open subset of $M \times M$ of full 
Riemannian measure. We will use the notation $d\mu(x)$ for the density of this measure.
Denote by $B(\widetilde x,T)$ the ball of radius $T$ in 
the universal cover of $M$ around a lift $\widetilde x$ of $x$.
By Proposition 4.3 in \cite{BeBo},
\begin{equation}  \label{ber_bot_eq}
 \int_M n_T(x,y) \,d\mu(y) \geq \hbox{\rm Vol}\, B(\widetilde x,T)
\end{equation}
for any $x \in M$.

Let $h=h(M)$ be the topological entropy of the geodesic flow. Manning showed in \cite{Manning} that
 $$
 h \geq \lim_{T \to \infty} \frac1T \log \hbox{\rm Vol}\, B(\widetilde x,T);
 $$
 the limit on the right hand side exists and is the same for all points 
$\widetilde x$ in the universal cover.
Ma\~n\'e improved this estimate in \cite{Mane} by showing that
\begin{equation}   \label{mane_eq}
h=\lim_{T\to\infty}\frac{1}{T}\log\int_{M\times M}n_T(x,y)\,d\mu(x)d\mu(y)
\end{equation}  
for any compact Riemannian manifold.

A Riemannian manifold has no conjugate points
if the exponential map at any point is a local diffeomorphism.
An equivalent characterization is that any two points in the universal 
cover are joined by exactly one geodesic.
The preceding relations are even stronger for compact manifolds with no conjugate points.  
In this case $n_T(x,y)$ is always finite, and for any $x,y \in M$ we have
\begin{equation} \label{rate}
h=\lim_{T \to \infty} \frac1T \log n_T(x,y).
\end{equation}
See for example Corollary 1.2 in \cite{Mane}. 

Equation~\eqref{rate} fails, in general, if $M$ has conjugate points.
There are examples where the growth of $n_T(x,x)$ is arbitrarily large 
for some exceptional points $x\in M$ \cite{BP96}. 
There are also examples in which
the limit in equation~\eqref{rate} is smaller than the topological entropy 
for an open set of configurations \cite{BP97}.

\subsection{Security and insecurity}\label{SSsec}
A set $B$ is a {\em blocking set} for a collection of geodesics if every geodesic 
belonging to the collection  passes through a point of $B$.

\begin{defin}  \label{secur_def} 
{\rm A configuration $(x,y)$ is  
{\em secure} if the collection $\Gamma(x,y)$ of connecting geodesics 
has a finite blocking set.
Otherwise the configuration is {\em insecure}. 
}
\end{defin}

A blocking set for $\Gamma(x,y)$ must lie in $M \setminus \{x,y\}$. It is also a
blocking set for $G(x,y)$.

If the 
configuration $(x,y)$ is secure, its
{\em security threshold}  $\sigma(x,y)$
is the minimal number of points needed to block $\Gamma(x,y)$.

\begin{defin}  \label{securmfld_def} 
{\rm  A Riemannian manifold is {\em secure} if
every configuration  is secure, {\em insecure} if some configurations are
insecure,  and {\em totally insecure} if every
configuration is insecure. It is {\em uniformly secure} if 
there is $s<\infty$ such that $\sigma(x,y)\le s$ for all $x,y\in M$.
The smallest such $s$ is
the {\em security threshold} of the manifold. 
}
\end{defin}

In a geometric optics interpretation, a configuration is secure
if one of the points can be shaded from the light emanating from the other
by a finite number of point screens. 
Another obvious interpretation of Definition~\ref{secur_def} suggested the 
name ``security": all geodesic paths connecting elements in a secure configuration can be 
monitored from a finite number of observation spots.
Since the security of configurations concerns the global properties of geodesics,
it is instructive to compare spaces from this viewpoint. 

\medskip

The framework of security arose in the context of polygonal billiards.
(The geodesics are the billiard orbits.) A polygon $P\subset\RR$
is {\em rational} if the angles of $P$ are commensurable with $\pi$.
The billiard orbits in a rational polygon $P$ correspond to the geodesics
in the associated {\em translation surface} $S=S(P)$. These surfaces
are endowed with flat metrics with special cone singularities.
Definitions~\ref{secur_def} and~\ref{securmfld_def} extend to polygons
and  translation surfaces. 

The results on security in this context are limited up to now to  
{\em lattice polygons} and {\em lattice translation surfaces}.
The lattice condition was introduced by Veech; the class of lattice 
translation surfaces partitions into two subclasses: arithmetic and
nonarithmetic \cite{GJ}. From the security standpoint, these subclasses
complement each other: An arithmetic translation surface is uniformly
secure; almost all configurations in a nonarithmetic translation surface
are insecure \cite{Gut05}. This dichotomy holds for billiard orbits in 
lattice polygons as well \cite{Gut05}.

\medskip

Gutkin and Schroeder investigated the security of 
compact, locally symmetric spaces in \cite{GS06}. 
Let $M$ be such a space. Then $M=S/\Ga$, where
$S$ is a simply connected symmetric space, and $\Ga$ is a discrete, cocompact group
of isometries freely acting on $S$. The space $S$  uniquely decomposes 
into a product  $S=S_0\times S_-\times S_+$  of 
simply connected symmetric spaces of {\em euclidean type, noncompact type, and 
compact type} respectively. If $M=S/\Ga$ where $S$ belongs to one of
the three types, we say that $M$ is a compact, locally symmetric space of that type.

\medskip

We summarize the results of \cite{GS06} in the following proposition.

\begin{prop}  \label{gut_schr_prop}
Let $M$ be a compact, locally symmetric space. Then the following statements hold:

\noindent a) If $M$ is of noncompact type then it is totally insecure.

\noindent b) If $M$ is of compact type then it has an open, dense, full measure
set of secure  configurations; it always has  insecure  configurations as well.

\noindent c) If $M$ is of euclidean type then 
it is uniformly secure, and its security threshold is bounded in terms of $\dim(M)$. 

\noindent d) If $M$ is secure, then it is of euclidean type.
\end{prop}
%

\section{Key Lemmas} \label{keysect}
The main results of this paper, which are presented in the next section, 
are based on the following lemmas. They have been obtained independently 
by Lafont and Schmidt \cite{LS}.

\begin{lem} \label{m&n}
Let $M$ be a compact Riemannian manifold, and let $\de$ be its injectivity radius.
Then for any $x,y\in M$ we have 
 $$
 m_T(x,y) \leq n_T(x,y) \leq \left(\frac{T}{2\delta}\right)^2 m_T(x,y).
 $$
\end{lem}
\begin{proof} The left hand inequality is obvious because $\Gamma_T(x,y) \subset G_T(x,y)$. 
We prove the other one by constructing a surjective map $\lambda:G_T(x,y)\to\Gamma_T(x,y)$ that is 
at most $(T/2\delta)^2:1$.

Observe firstly that if $t$ and $t'$ are consecutive times when a geodesic passes 
through a point, $z$, 
then $2\delta\le|t - t'|$. Indeed, the geodesic must travel (at least) from the center 
of the ball $B(z,\delta)$ to its boundary, and then back. 
An immediate corollary is that the set of times when a geodesic passes through $z$ is discrete.

Let $\gamma: [a,b] \to M$ be a geodesic that belongs to $G_T(x,y)$. This means that
$\gamma(a) = x$ and $\gamma(b) = y$. 
Let $b'$ be the leftmost element of $(a,b]$ such that $\gamma(b') = y$.
Let $a'$ be the rightmost element of $[a,b')$ such that $\gamma(a') = x$. 
Observe that the restriction of $\gamma$ to $[a',b']$ belongs to $\Gamma_T(x,y)$ 
and that $[a',b']$ is the leftmost subinterval of $[a,b]$ with this property. 
Furthermore, if $\gamma$ already belongs to $\Gamma_T(x,y)$, we will have $a' = a$ and $b' = b$.

Define $\lambda(\gamma)$ to be the restriction of $\gamma$ to $[a',b']$. The resulting map
$\lambda:G_T(x,y)\to\Gamma_T(x,y)$ is surjective because it is the identity 
on $\Gamma_T(x,y)$. 

Let $\gamma': [a',b'] \to M$  belong to $\Gamma_T(x,y)$. We will investigate the preimage
$\la^{-1}(\gamma')\subset G_T(x,y)$. Let $\gamma:{ \Bbb R} \to M$ be the 
extension of $\gamma'$ to a complete geodesic. The geodesics in $\la^{-1}(\gamma')$
are the restrictions of $\gamma$ to intervals of the form $[a,b]$ 
such that $\gamma(a) = x$ and $\gamma(b) = y$. We must also have 
$a \leq a' < b' \leq b$, and $|b - a| \leq T$.
These inequalities ensure that $a \in [b' - T,a']$ and $b \in [b',a' + T]$; both of these intervals have length $\le T$. 
We see from the remarks at the beginning of the proof that there can be at most $T/2\delta$
possible choices for $a \in [b' - T,a']$ with $\gamma(a) = x$. Similarly there can be at most $T/2\delta$
possible choices for $b \in [b',a' + T]$ with $\gamma(b) = y$. It follows that $\gamma'$
has at most $(T/2\delta)^2$ preimages under $\la$, as desired.
\end{proof}


\begin{lem} \label{basiclem}
Let $M$ be a compact Riemannian manifold.
Let $\Gamma_T(x,y;z)$ denote the subset of geodesics in $\Gamma_T(x,y)$ that pass 
through  $z \neq x,y$. Let $m_T(x,y;z)$ be the number of geodesics in $\Gamma_T(x,y;z)$.
 Then 
 $$
 m_T(x,y;z) \leq m_{T/2}(x,z) + m_{T/2}(z,y).
 $$
\end{lem} 
\begin{proof} We construct a 1:1 map $\mu:\Gamma_T(x,y;z)\to\Gamma_{T/2}(x,z) 
\cup \Gamma_{T/2}(z,y)$.

Suppose $\gamma:[a,b] \to M$ is a geodesic in $\Gamma_T(x,y;z)$. Let $\gamma_1$ be the initial segment of 
$\gamma$ that ends at the first time when $\gamma$ passes through $z$; let 
$\gamma_2$ be the terminal segment of $\gamma$ that 
starts at the last time when $\gamma$ passes through $z$. 
Observe that  $\gamma_1$ and $\gamma_2$ both belong to $\Gamma(x,z) \cup \Gamma(z,y)$,
and their domains are subintervals of $[a,b]$
that intersect in at most one point. 

Since the total length of $\gamma_1$ and $\gamma_2$ is at most $T$, 
at least one of them has length $\leq T/2$. Set $\mu(\ga)=\ga_1$ if 
$\ga_1$ has length $\leq T/2$ and set $\mu(\ga)=\ga_2$ otherwise. We will now show 
that $\mu$ is 1:1.

A geodesic in $\Gamma(x,z)$ can be extended past $z$ to obtain a geodesic 
in $\Gamma(x,y)$ in at most one way: we must stop at the first time (if there is one) 
when the geodesic reaches $y$. 
Similarly, a geodesic in $\Gamma(z,y)$ can be extended past $z$ 
in at most one way to obtain a geodesic in  $\Gamma(x,y)$. 

Thus each $\ga_1\in\Gamma_{T/2}(x,z)$ 
can be the intial segment for at most one $\ga\in\Gamma(x,y)$. 
Similary, each $\ga_2\in\Gamma_{T/2}(z,y)$) can be the terminal segment 
for at most one $\ga\in\Gamma(x,y)$.
Therefore the mapping  
$\mu:\Gamma_T(x,y;z)\to\Gamma_{T/2}(x,z) \cup \Gamma_{T/2}(z,y)$
is injective.
\end{proof}

\section{Main Results}   \label{mainsect}
Our main technical result relates the security threshold of a manifold with 
the growth rate of joining geodesics.
\begin{prop}   \label{main_prop}
Let $M$ be a compact Riemannian manifold, and let $\de$ be its injectivity radius.
Suppose that $M$ is uniformly secure, and let $s=s(M)$ be the security threshold.
Then for any $x, y \in M$ we have 
$$
  n_T(x,y) < \frac{s}2 \cdot \left( \frac T\delta \right)^{3+\log_2s}.
$$
\begin{proof}
By the uniform security assumption, any configuration $(p,q)\in M \times M$ has a 
blocking set 
$B(p,q)\subset M\setminus\{p,q\}$ containing at most $s$ points.
{From} Lemma~\ref{basiclem}, we have
\begin{equation}   \label{reduc_eq}
m_T(x,y) \leq \sum_{z \in B(x,y)}\left( m_{T/2}(x,z) + m_{T/2}(z,y) \right).
\end{equation}
The right hand side of this inequality has the 
form $\sum m_{T/2}(p_i,q_i)$, where the
$(p_i,q_i)$ run through a finite set, $\pp_1$, of pairs of points in $M$ with 
at most $2s$ elements.

Iterating the estimate above $k$ times, we obtain a sequence, $\pp_k$, 
of subsets of $M \times M$, such that $\pp_k$ has at most $(2s)^k$ elements and
\begin{equation}   \label{induc_eq}
m_T(x,y) \leq \sum_{(p_i,q_i)\in\pp_k} m_{T/2^k}(p_i,q_i).
\end{equation}

Observe that for any $L>0$ we have $m_L(p_i,q_i)\le\sup_{p,q\in M}m_L(p,q)$.
Combining this with inequality~\eqref{induc_eq} yields
\begin{equation}   \label{estim_eq}
m_T(x,y) \leq (2s)^k \sup_{p,q\in M}m_{T/2^k}(p,q).
\end{equation}

Now choose $k\in\N$ so that $\log_2(\frac{T}{\de})< k \le 1+\log_2(\frac{T}{\de})$, or equivalently
$$
\frac{\de}2 \le \frac{T}{2^k} < \de.
$$
Since $T/2^k < \delta$, we have $m_{T/2^k}(p,q)\le 1$ for any $p,q\in M$. 
Substituting this into inequality~\eqref{estim_eq},
we obtain
\begin{equation}   \label{final_eq}
m_T(x,y) < 2s \left(\frac{T}{\de}\right)^{(1+\log_2s)}.
\end{equation}
Combining this inequality with Lemma~\ref{m&n} yields the claim.
\end{proof}
\end{prop}

Proposition~\ref{main_prop} directly implies that a uniformly secure manifold has  
uniform polynomial growth of geodesics. The following theorem states this precisely.
Recall that a positive function $f(T)$ has polynomial growth if 
there are constants $C$ and $d$ such that $f(T) \leq CT^d$ as $T \to \infty$. 
The minimum possible $d$ is the degree of the polynomial growth.

\begin{thm}   \label{polin_thm}
Let $M$ be a compact Riemannian manifold. If $M$ is uniformly secure, 
then there are positive constants
$C$ and  $d$ such that for any pair $x,y\in M$ we have
$$
n_T(x,y) \le C T^d.
$$
\end{thm}

This inequality yields strong consequences for the geometry and topology of secure manifolds.

\begin{thm}  \label{entrop_thm}
Let $M$ be a compact Riemannian manifold that is uniformly secure. 
Then the  topological entropy of the geodesic flow for $M$ is zero, 
and the fundamental group of $M$ is virtually nilpotent. 

If, in addition, $M$ has no conjugate points, then $M$ is flat.
\end{thm}
\begin{proof}

Our first claim is immediate from equation~\eqref{mane_eq} 
and Theorem~\ref{polin_thm}. Set
$$
N(x,T)=\int_M n_T(x,y) \,dy.
$$
By Theorem~\ref{polin_thm}, $N(x,T)$ has   
polynomial growth as $T \to \infty$ for any $x \in M$. 
By equation~\eqref{ber_bot_eq}, this implies polynomial growth 
as $T \to \infty$ of the volume of the ball
$B(\widetilde x,T)\subset\widetilde M$ of radius $T$ in the 
universal cover of $M$, which in turn
immediately implies polynomial growth of the fundamental group of $M$. 
The latter means that for a finite set $S$ of generators for $\pi_1(M)$ the number of
elements of $\pi_1(M)$ expressible as words of length $\leq n$ in the alphabet $S$
has polynomial growth. Polynomial growth for one finite set of generators implies 
polynomial growth of the same degree for any finite set of generators.

By a famous theorem of Gromov, a finitely generated group with polynomial growth 
is virtually nilpotent \cite{Grom}, i.e.~it has a finite index nilpotent subgroup. 
This proves our
second claim. Our last claim now follows from a recent result of 
Lebedeva \cite{Leb}: if the fundamental group of
a compact Riemannian manifold with no conjugate points has polynomial growth,
then the manifold is flat.
\end{proof}

%

Theorem~\ref{entrop_thm} has an immediate corollary.
\begin{corol}   \label{admit_cor}
Let $M$ be a compact manifold. If $\pi_1(M)$
is not virtually nilpotent, then $M$ does not admit a  
uniformly secure Riemannian metric.
\end{corol}

Our next result concerns manifolds with no conjugate points. 
It has been obtained independently by Lafont and Schmidt \cite{LS}.


\begin{thm}   \label{nocp_thm}
Let $M$ be a compact Riemannian manifold with no conjugate points 
whose geodesic flow has positive topological entropy. Then $M$ is 
totally insecure.
\end{thm}
\begin{proof} 
Let $p,q \in M$ be arbitrary. From equation~\eqref{rate} and Lemma~\ref{m&n}, we have
\begin{equation} \label{con_entr}
\lim_{T\to \infty} \frac1T \log m_T(p,q) = h.
\end{equation}
This is where we use the no conjugate points property.

Suppose now that the claim is false. Then there are two points $x,y \in M$ 
such that any geodesic from $x$ to $y$ is blocked by one of a finite number of points, say $z_1,\dots,z_n$.  Lemma~\ref{basiclem} tells us that
 \begin{equation}\label {T&T/2}
 m_T(x,y) = \sum_{i = 1}^n \left( m_{T/2}(x,z_i) +  m_{T/2}(z_i,y)\right)
 \end{equation}
 for any $T > 0$. Consider a number $\alpha > 1$. By \eqref{con_entr}, there is a
$T_0 = T_0(\alpha)$ such that for any $T \geq T_0$ and any $p,q \in \{x,y,z_1,\dots,z_n\}$ 
we have
 $$
 h/\alpha < \frac1T \log m_T(p,q) < \alpha h,
 $$
 which is equivalent to
 $$
 e^{hT/\alpha} < m_T(p,q) < e^{hT\alpha} .
 $$
 Combining this with equation~\eqref{T&T/2}, we obtain
 $$
 e^{hT/\alpha} < 2n e^{hT\alpha/2}
 $$
 for all large enough $T$. Now choose $\alpha < \sqrt{2}$. Then the exponent on the left is 
larger than the exponent on the right, and the inequality is absurd for large $T$.
\end{proof}

%
%
%
%

%
%

\medskip

Theorem~\ref{nocp_thm} and the variational principle for entropy \cite{Goodw} 
immediately imply the following proposition.

\begin{corol}   \label{pos_ent_cor}
Let $M$ be a compact Riemannian manifold with no conjugate points. 
If the geodesic flow of $M$ has positive metric entropy,
then the manifold $M$ is totally insecure.
\end{corol}

\medskip

Any compact, locally symmetric space of noncompact type satisfies the hypotheses of 
Theorem~\ref{nocp_thm}. This gives a new proof of case (a) in 
Proposition~\ref{gut_schr_prop}. The original proof in \cite{GS06} 
used an entirely different approach.

For manifolds of nonpositive curvature we obtain a dichotomy along
the lines security/insecurity.

\begin{corol}     \label{neg_curv_cor}
A compact Riemannian manifold of nonpositive curvature is either totally
insecure or uniformly secure. In  the latter case the manifold is flat; 
its security threshold is bounded above in terms of the dimension
of manifold. 
\end{corol}
\begin{proof}
Let $M$ be a compact Riemannian manifold of nonpositive curvature.
By Pesin's formula (Corollary 3 in \cite{Pes}),
the following dichotomy holds: either i)
the geodesic flow for $M$ has positive entropy, or ii) $M$ is flat.

If i) holds, then $M$ satisfies the hypotheses of Corollary~\ref{pos_ent_cor};
thus, $M$ is totally insecure. In case ii), $M$ is 
a compact, locally symmetric space of euclidean type. The claim follows from 
case (c) of Proposition~\ref{gut_schr_prop}.
\end{proof}

\section{Examples, conjectures and open problems}    \label{conj&probs}  
We begin by producing examples of insecure manifolds with conjugate points. 

\begin{lem} \label{prod_lem}
A configuration in a Riemannian product $M \times M'$ is secure if and
only if its projections to $M$ and $M'$ are secure.
\end{lem}

\begin{proof}
Geodesics in $M \times M'$ are products of geodesics in $M$ and geodesics in $M'$.
More precisely, if $\gamma:[a,b] \to M$ and
$\gamma':[a',b'] \to M'$ are unit speed geodesics, then the unit speed geodesic 
$\gamma \times\gamma'$ in $M$ of length
$L = \sqrt{(b-a)^2 + (b'-a')^2}$ is defined by
 $$
 (\gamma \times \gamma')(t) = (\gamma(a+(b-a)t/L),\gamma'(a'+(b'-a')t/L)
 $$
for $0 \leq t \leq L$. All geodesics in $M \times M'$ arise in this way.

Let $\xi = (x,x')$ and $\eta = (y,y')$ be two points in $M\times M'$. The set
$G(\xi,\eta)$ consists of the products of geodesics in $G(x,y)$ and
$G(x',y')$. Let $\Gamma_0(\xi,\eta)$ be the subset of $G(\xi,\eta)$ formed by the 
products of geodesics in $\Gamma(x,y)$ and $\Gamma(x',y')$.
Then $\Gamma_0(\xi,\eta)\subset\Gamma(\xi,\eta)$. (There may
be more geodesics in $\Gamma(\xi,\eta)$, for example products of
geodesics in $\Gamma(x,y)$ and in $G(x',y') \setminus \Gamma(x',y')$.)

Suppose now that $x$ and $y$  form a secure configuration in $M$ with
blocking set $B \subset M \setminus \{x,y\}$,   and that $x'$ and $y'$
form a secure configuration in $M'$ with blocking set 
$B' \subset M' \setminus \{x',y'\}$. Then $B$ and $B'$
also block all the geodesics in $G(x,y)$ and $G(x',y')$ respectively.
If $x \neq y$ and $x' \neq y'$, the set $B \times B'$ will block all
geodesics in $G(\xi,\eta)$. If $x = y$, we will also use 
the blocking points
$(x,b')$, where $b' \in B'$, in order to block geodesics 
from $\xi$ to $\eta$
that have constant projection to the first factor. Similarly if $x' = y'$, we
include the points $(b,x')$, where $b \in B$, in our blocking set.
In all cases, we produce a finite set $Z\subset M\times M'$, such that
$B \times B' \subset Z\subset M\times M' \setminus \{\xi,\eta\}$,
 blocking all geodesics in  $G(\xi,\eta)$. 

Conversely, if $\xi$ and $\eta$ form a secure configuration, there is a finite set of
points $(z_1,z'_1), \dots, (z_n,z'_n)$ in $M \times M'$ which block
all geodesics in $\Gamma_0(\xi,\eta)$. Then $\{z_1,\dots, z_n\}$ is a blocking set for
$\Gamma(x,y)$ and $\{z'_1,\dots, z'_n\}$ is a blocking set
for $\Gamma(x',y')$.
\end{proof}

The following criteria for security and insecurity follow immediately from the preceding lemma 
and its proof.

\begin{prop}    \label{riem_prod_prop}
Let $M$ and $M'$ be Riemannian manifolds.
Then the following statements hold.

\noindent 1) If one of $M$ and $M'$ is (totally) insecure, then $M \times M'$
is  (totally) insecure.

\noindent 2. $M \times M'$ is  (uniformly) secure if and only if both $M$ and $M'$ are
(uniformly) secure.
\end{prop}

\begin{rem}   \label{cover_rem}
{\em Lemma~\ref{prod_lem} and
Proposition~\ref{riem_prod_prop} are closely related to 
Proposition 5 and Corollary 4 in \cite{GS06}. They extend to fibre bundles such that 
each point of the base has
a neighborhood, $O$, over which the fibration is the Riemannian product of $O$ and the fibre.
We leave the details to the reader.  
}
\end{rem}  

\medskip

Combining Proposition~\ref{riem_prod_prop} with Theorem~\ref{nocp_thm} immediately gives the 
following result.

\begin{corol}  \label{nocp_corol}
Let $M$ be a compact Riemannian manifold with no conjugate points
whose geodesic flow has positive topological entropy. Let $M'$ be a 
Riemannian manifold. Then the Riemannian product $M \times M'$ is totally insecure.
\end{corol}
This proposition provides us with examples 
of totally insecure manifolds that have conjugate points.

The hypothesis of no conjugate points in
Theorem~\ref{nocp_thm} is used only to ensure that equation~\eqref{rate} holds.
One can strengthen equation~\eqref{rate} by showing 
that the convergence to the limit is uniform for all $p,q \in M$. The proof of 
Theorem~\ref{nocp_thm} does not require this uniformity. 

In fact, in order to prove Theorem~\ref{nocp_thm}, we need considerably 
less than equation~\eqref{rate}. 
It suffices to have positive constants $h_1$ and $h_2$ such that $h_2 < 2h_1$,
and such that for arbitrary pairs $p,q\in M$ and all large enough $T$ we have
the inequality
\begin{equation}    \label{up_low_eq}
e^{h_1T} \leq n_T(p,q) \leq  e^{h_2T}.
\end{equation}    
Suppose that inequality~\eqref{up_low_eq} holds. Then  
all configurations in $M$ are insecure, whether there are conjugate points or not.  
We believe that many Riemannian manifolds satisfy this condition,
but the only examples that we know presently are the  manifolds with no conjugate points. 

\medskip

Following another approach, we can replace  inequality~\eqref{up_low_eq} by
\begin{equation}    \label{limits_eq}
\sup_{x,y \in M}\limsup_{T \to \infty}\frac1T\log n_T(x,y) 
                         <  2\liminf_{T \to \infty}\frac1T\log n_T(p,q) < \infty.
\end{equation}    
Similarly, inequality~\eqref{limits_eq} implies that the pair $p,q\in M$ is insecure.  
We believe that inequality~\eqref{limits_eq} frequently holds, 
but do not know any examples outside the realm of manifolds with no conjugate points.

\medskip

Examples show that one cannot simply discard the hypothesis of no conjugate points from 
Theorem~\ref{nocp_thm}. In fact,  the claim does not extend to all manifolds 
whose geodesic flows have positive topological entropy. There are  compact Riemannian manifolds
with positive topological entropy that contain secure configurations. 

We outline  such an example; it was used for a different purpose in \cite{BP97}. 
It is a special case of a construction due to Weinstein. See \cite{Be3},
Appendix C. Choose a Riemannian metric $\widehat g$
on $S^2$ whose geodesic flow has positive topological entropy. Such metrics exist; the 
first was given by Knieper and Weiss \cite{KW}. Now choose a smooth family of metrics $g_t$,
$-\pi/2 \leq t \leq \pi/2$, such that $g_t = \widehat g$ when $|t| \leq \pi/6$, and $g_t$ is
the standard metric on the unit sphere in $\Bbb R^3$ when $\pi/3 \leq |t| \leq \pi/2$. 
Choose also a scaling function $\rho: [-\pi/2,\pi/2] \to \Bbb R$ 
such that $\rho(t) = 1$ when $|t| \leq \pi/6$, $\rho(t) > 0$ when $\pi/6 \leq |t| \leq \pi/3$,
and $\rho(t) = \cos t$ when $\pi/3 \leq |t| \leq \pi/2$.

Define a metric on $[-\pi/2,\pi/2] \times S^2$ that is degenerate on the boundary as follows. 
Let $t \in [-\pi/2,\pi/2]$ denote the first coordinate. We require:
\begin{enumerate}
\item $\|\partial/\partial t\| \equiv 1$;
\item $\partial/\partial t$ is everywhere orthogonal to $\{t\} \times S^2$;
\item the restriction of the metric to the vectors tangent to 
$\{t\} \times S^2$ is $\rho(t)g_t$.
\end{enumerate}
The result is a smooth metric $g$ on $S^3$. The spheres $\{\pm \pi/2\} \times S^2$ collapse to 
points, $p$ and $p'$, which we call poles. The metric $g$ coincides with the 
standard round metric on $S^3$ near the poles.

The spheres $\{t\} \times S^2$ for $|t| \leq \pi/6$ are totally geodesically isometrically embedded copies of $(S^2,\widehat g)$; 
this ensures that the geodesic flow for the metric $g$ on $S^3$ has positive entropy. 
The geodesics which pass through the  poles are circles of length $2\pi$ 
and are orthogonal to the ``spheres of latitude''. 
They are meridians just as in the usual round metric on $S^3$; 
the only difference is that the spheres of latitude have an unusual 
geometry away from the poles. All meridians that leave one pole focus at the 
opposite pole after distance $\pi$. 
The pairs $(p,p)$ and $(p',p')$ are secure configurations because their connecting  geodesics 
are blocked by the opposite pole.

\medskip

%
We conclude with some conjectures related to the results of this paper.

Presently, the only 
examples of compact Riemannian manifolds that are secure are the
flat manifolds. We believe that this is not an accident and Theorems~\ref{polin_thm} and 
\ref{entrop_thm} are evidence for:

\begin{conj}   \label{secur_flat_conj}
A  compact Riemannian manifold is secure if and only if it is  flat.
\end{conj}

Our next conjecture concerns classical Riemannian geometry but 
is motivated by security questions, as we explain below.

\begin{conj}    \label{noconju_flat_conj}
Let $M$ be a  compact Riemannian manifold with no conjugate points.
Then either $M$ is flat or its geodesic flow has positive topological entropy.
\end{conj}

If Conjecture~\ref{noconju_flat_conj} is true, then the  alternative for 
manifolds with nonpositive curvature
stated in Corollary~\ref{neg_curv_cor} would extend to compact Riemannian manifolds
with no conjugate points. More precisely,  the following dichotomy would hold: 
A compact manifold with no conjugate points is either

\noindent i) flat, uniformly secure, 
and has zero topological entropy for its geodesic flow; or

\noindent ii)  
not flat, totally insecure,
and has positive topological entropy for its geodesic flow. 

\medskip

Finally, we suspect  that there are 
many examples of totally insecure manifolds other than those described in this paper.  
In fact, we believe that ``most'' Riemannian manifolds are totally insecure.

More conjectures related to security in Riemanian manifolds can be found in \cite{LS}.


\begin{thebibliography}{99}
%


\bibitem{BeBo} M. Berger and R. Bott, {\em Sur les vari\'et\'es \`a courbure strictement positive}, 
Topology {\bf 1} (1962), 302 -- 311.

\bibitem{Be3} A.L. Besse, {\em Manifolds all of whose geodesics are closed}, Springer-Verlag, New York 1978.





\bibitem{BP96} K. Burns and G. Paternain, {\em On the growth 
of the number of geodesics joining two points}, 
Pitman Res. Notes Math. {\bf 362} (1996), 7 -- 20.



\bibitem{BP97} K. Burns and G. Paternain, {\em  
Counting geodesics on a Riemannian manifold and topological entropy of geodesic flows}, 
Erg. Theo. \& Dyn. Sys. {\bf 17} (1997), 1043 -- 1059.





















\bibitem{Grom} M. Gromov, {\em Groups of polynomial growth and expanding maps},
Inst. Hautes \'Etudes Sci. Publ. Math. {\bf 53} (1981), 53 -- 73.

\bibitem{Goodw} L.W. Goodwyn,  {\em Topological entropy bounds measure-theoretic entropy},  Proc. A. M. S.  {\bf 23}  (1969), 679 -- 688.









\bibitem{Gut05} E. Gutkin, {\em Insecure configurations in lattice translation surfaces, 
with applications to polygonal billiards}, Discr. Cont. Dyn. Sys. {\bf A 16} (2006), 367 -- 382.


\bibitem {GJ} E. Gutkin and C. Judge, {\em Affine mappings of translation surfaces: 
geometry and  arithmetic}, Duke Math. J. {\bf 103} (2000), 191 -- 213.  


\bibitem{GS06} E. Gutkin and V. Schr\"oder, {\em Connecting geodesics and 
security of configurations in compact 
locally symmetric spaces}, Geometriae Dedicata {\bf 118} (2006), 185 -- 208.






 

\bibitem{KW} G. Knieper and H. Weiss, {\em A surface with positive curvature 
and positive topological entropy},  J. Diff. Geom.  {\bf 39}  (1994), 229 -- 249.




\bibitem{LS} J.-F. Lafont and B. Schmidt, {\em Blocking light in compact Riemannian manifolds},
preprint (2006).

\bibitem{Leb} 
N.D. Lebedeva, {\em  On spaces of polynomial growth with no conjugate points},
St. Petersburg Math. J.  {\bf 16}  (2005), 341 -- 348.



\bibitem{Mane} R. Ma\~n\'e, {\em On the topological entropy of geodesic flows}, 
J. Diff. Geom.  {\bf 45} (1997), 74 -- 93. 

\bibitem{Manning} A. Manning, {\em Topological entropy for geodesic flows}, Ann. Math.
{\bf 110}  (1979), 567 -- 573.











\bibitem{Pes} Ya. Pesin, {\em Formulas for the entropy of the geodesic flow on 
a compact Riemannian manifold without conjugate points}, Math. Notes {\bf 24} (1978), 796 -- 805.


























































\end{thebibliography}
\end{document}